\documentclass [a4paper,twoside,12pt]{article}
\usepackage[utf8]{inputenc}
\usepackage{amsmath,amsfonts,amssymb}
\usepackage{vmargin,graphicx,theorem}
\usepackage[english]{babel}
\usepackage{enumerate}
\usepackage{color}
\usepackage{pstricks}
\usepackage{mdwlist}
\setpapersize[portrait]{A4}
\setmarginsrb{1.5cm}{1cm}{1.5cm}{2.5cm}{1.5cm}{0cm}{0.5cm}{2cm}
\newcommand{\rr}{\ensuremath{\mathbb{R}}}

\newcommand{\dR}{\ensuremath{\mathbb{R}}}

\newtheorem{ethm}{Theorem}[section]

\newtheorem{eprop}[ethm]{Proposition}

\newtheorem{elem}[ethm]{Lemma}

\newtheorem{erem}[ethm]{Remark}

\newcommand{\proofend}{~$\rhd$}
\newcommand{\proofbegin}{~$\lhd$}

\newenvironment{eproof}
               {\noindent {\emph{\textbf{Proof}}}\\\proofbegin~}
               {\proofend\\}

\renewcommand{\phi}{\varphi}

\renewcommand{\geq}{\geqslant}

\newcommand{\R}{\dR}

\newcommand{\beq}{\begin{equation}}\newcommand{\eeq}{\end{equation}}
\begin{document}

\title{Uniform convergence to equilibrium for granular media}
\author{ Fran{\c c}ois Bolley\thanks{Ceremade, Umr Cnrs 7534, Universit\'e Paris-Dauphine, Place de Lattre de Tassigny, F-75775 Paris cedex 16. bolley@ceremade.dauphine.fr}, Ivan Gentil\thanks{Institut Camille Jordan, Umr Cnrs 5208, Universit\'e Claude Bernard Lyon 1, 43 boulevard du 11 novembre 1918, F-69622 Villeurbanne cedex. gentil@math.univ-lyon1.fr}\, and Arnaud Guillin\thanks{Institut Universitaire de France and Laboratoire de Math\'ematiques, Umr Cnrs 6620, Universit\'e Blaise Pascal, Avenue des Landais, F-63177 Aubi\`ere cedex. guillin@math.univ-bpclermont.fr}}

\date{\today}

\maketitle

\abstract{We study the long time asymptotics of a nonlinear, nonlocal equation used in the modelling of granular media. We prove a uniform exponential convergence to equilibrium for degenerately convex and non convex interaction or confinement potentials,  improving in particular results by J.~A.~Carrillo, R.~J.~McCann and C.~Villani. The method is based  on studying the dissipation of the Wasserstein distance between a solution and the steady state. }
\bigskip

\bigskip

\section*{Introduction}

We consider the problem of convergence to equilibrium for the nonlinear equation
 \begin{equation}
 \partial_t\mu_t=\Delta \mu_t+\nabla \cdot (\mu_t(\nabla V+\nabla W*\mu_t)) \qquad\qquad t>0,x\in\R^n.
 \label{eq-gran}
 \end{equation}
 This equation preserves mass and positivity and we shall be concerned with solutions which are probability measures on $\rr^n$ at all times. It is used in the modelling of space-homogeneous granular media (see~\cite{BCCP}), where it governs the evolution of the velocity distribution $\mu_t (x)$  of a particle under the effects of diffusion, a possible exterior potential $V$ and a mean field interaction through the potential $W$; we shall keep the variable $x$ instead of $v$ (for the velocities) for notational convenience.
 
Steady states may exist as a result of a balance between these three effects, and we are concerned with deriving rates of convergence of solutions towards them. Following~\cite{BCCP}, this issue has raised much attention in the last years and has been tackled by a particle approximation and logarithmic Sobolev inequalities in~\cite{malrieu03}, by an entropy dissipation method in~\cite{cmcv-03,cordero-gangbo-houdre}  and by contraction properties in Wasserstein distance in~\cite{cmcv-06,cgm-08,bgm} (see also \cite{CT03, calvez-carrillo10} for related works in dimension one).  The entropy method is based on studying the time derivatives of a Lyapunov function $F$ of the equation (called entropy or energy), on the interpretation due to F. Otto of~\eqref{eq-gran} as a gradient flow of $F$  (see~\cite{cmcv-06,ambrosio-gigli-savare}) and on the notion of convexity for $F$ due to R.~J. McCann (see~\cite{mccann-advances}). 

\medskip

When $V$ and $W$ are uniformly convex, solutions converge exponentially fast to equilibrium, but the case of interest of~\cite{BCCP} is $V=0$ and $W(x)=|x|^3$, whose convexity degenerates at~$0$. For this case, only a polynomial rate, or exponential but depending on the initial data, was obtained in~\cite{cmcv-03,cmcv-06,cgm-08}. In the present paper we prove a {\it uniform exponential} convergence in Wasserstein distance of all solutions to the steady state. The method, introduced in the linear case in~\cite{BGG11}, is based on comparing the Wasserstein distance with its dissipation along the evolution. 

\smallskip

In Section~\ref{sec1} we derive the dissipation of the  Wasserstein distance between solutions and easily deduce the classical contraction results.  Section~\ref{sec2} is devoted to cases when the convergence is driven by the interaction potential $W$, with or without exterior potential $V$: in particular we  prove the first result of uniform exponential convergence to equilibrium for degenerately convex interaction potentials and no exterior potential. In Section~\ref{sec3} we give conditions to get an exponential convergence to equilibrium with both potentials being non convex.


\section{Dissipation of the Wasserstein distance}\label{sec1}

Let $P_2(\rr^n)$ be the set of Borel probability measures on $\rr^n$ with $ \int _{\rr^n}  \vert x \vert^2 \, d\mu<\infty$. The Wasserstein distance between two measures $\mu$ and $\nu$ in $P_2(\rr^n)$ is defined as
$$
W_2( \mu, \nu) = \inf_{\pi} \Big( \iint_{\rr^{2n}} \vert x-y \vert^2 \, d\pi(x,y) \Big)^{1/2}
$$
where $\pi$ runs over the set of joint Borel probability measures on $\rr^{2n}$ with marginals $\mu$ and $\nu$. It defines a distance on $P_2(\rr^n)$ which metrizes the narrow convergence, up to a condition on moments. In the present work, convergence estimates will be given in terms of this distance, but interpolation estimates can turn such weak convergence estimates into strong convergence estimates. By the Brenier Theorem, if $\mu$ is absolutely continuous with respect to the Lebesgue measure, then there exists a convex function $\phi$ such that $\nabla \phi\#\mu=\nu$, that is, 
$
\int_{\rr^{n}} g \, d\nu=\int_{\rr^{n}} g(\nabla \phi) \,d\mu
$
for every bounded  function $g$; moreover
$$
W_2^2(\mu,\nu)=\int_{\rr^{n}} |\nabla\phi(x) - x|^2 \, d\mu(x)
$$
and $\nabla \phi^* \#\nu = \mu$ for the Legendre transform $\phi^*$ of $\phi$ if also $\nu$ is absolutely continuous with respect to the Lebesgue measure. We refer to~\cite{portoercole,ambrosio-gigli-savare,villani-book1} for instance for these~notions. 

\medskip

We shall assume that $V$ and $W$ are $\mathcal C^2$ potentials on $\rr^n$, respectively $\alpha$ and $\beta$-convex with $\alpha,\beta \in\rr$, in the sense that $\nabla^2 V(x) \geq \alpha$ and $\nabla^2 W (x) \geq \beta$ for all $x \in \rr^n$, as quadratic forms on $\rr^n$. Moreover we assume that the interaction potential $W$ is even and that both $V$ and $W$ satisfy the doubling condition
\begin{equation}
\label{eq-doubling}
V(x+y) \leq C (1 + V(x) + V(y))
\end{equation}
for all $x,y \in \rr^n$, and analogously for $W$.

We shall consider solutions which are gradient flows in the space $P_2(\rr^n)$ of the free energy 
\begin{equation}\label{eq:F}
 F(\mu) = \int_{\rr^n} \mu \, \log \mu \, dx + \int_{\rr^n} V \, d\mu + \frac{1}{2} \iint_{\rr^{2n}} W(x-y) \, d\mu(x) \, d\mu(y),
\end{equation}
 as developed as follows in~\cite{ambrosio-gigli-savare}: Let $\mu_0$ be an initial datum in $P_2(\rr^n)$. Then, by~\cite[Ths. 4.20 and 4.21]{daneri-savare} or~\cite[Th. 11.2.8]{ambrosio-gigli-savare}, there exists a unique curve $\mu = (\mu_t)_t \in \mathcal C([0, + \infty[, P_2(\rr^n)$), locally Lipschitz on $]0, + \infty[$, satisfying the evolution variational inequality
$$
\frac{1}{2} \frac{d}{dt} W_2^2(\mu_t, \sigma) \leq F(\sigma) - F(\mu_t) - \frac{\alpha + \min\{\beta,0\}}{2} \, W_2^2(\mu_t, \sigma)
$$
for almost every $t>0$ and all probability measure $\sigma$ in the domain of $F$. For all $t>0$ the solution $\mu_t$ has a density with respect to the Lebesgue measure. Moreover the curve $
\mu$ satisfies the continuity equation
$$
\partial_t \mu_t + \nabla \cdot (\mu_t \, v_t) = 0, \qquad\qquad t>0,x\in\R^n
$$ 
in the sense of distributions, where the velocity field $v_t$ satisfies
$$
-\mu_t v_t = \nabla \mu_t + \mu_t \, \nabla V + \mu_t \, (\nabla W * \mu_t).
$$
In other words $\mu$ is a solution to~\eqref{eq-gran}, and the curve $\mu = (\mu_t)_t$ will be called the solution with initial datum $\mu_0 \in P_2(\rr^n)$. Finally $t \mapsto \int \vert v_t \vert^2 \, d\mu_t \in L^{\infty}_{loc} ([0, + \infty[)$ so the curve $\mu : \, ]0, + \infty[ \to P_2(\rr^n)$ is absolutely continuous (see~\cite[Th. 8.3.1]{ambrosio-gigli-savare}); moreover, if  $\nu$ is another such solution with initial datum $\nu_0$ and associated velocity field $w_t$, then by~\cite[Th. 23.9]{villani-book1} or~\cite[Th. 8.4.7]{ambrosio-gigli-savare}
$$
\frac{1}{2} \frac{d}{dt} W_2^2(\mu_t, \nu_t) = - \int_{\rr^n} (\nabla \varphi_t(x) - x) \cdot v_t(x) \, d\mu_t(x) - \int_{\rr^n} (\nabla \varphi_t^*(x) - x) \cdot w_t(x) \, d\nu_t(x)
$$
for almost every $t>0$; here $\phi_t$ is a convex function on $\rr^n$ such that $\nabla \phi_t \#\mu_t =\nu_t$ and $\nabla \phi_t^*  \#\nu_t =\mu_t.$ Then one can perform a ``weak" integration by parts as in~\cite[Th. 1.5]{lisini} or~\cite[Lem.~13]{cmcv-06} to bound from above the right-hand side by
$$
- \int_{\rr^n} \Big( \Delta\phi_t (x)+\Delta\phi_t^*(\nabla\phi_t(x))-2n +({\mathcal A}(\nabla\phi_t(x),\nu_t)-{\mathcal A}(x,\mu_t)) \cdot (\nabla\phi_t(x)-x) \Big) \,d\mu_t(x).
$$
Here $\Delta\phi$ is the trace of the Hessian of a convex map $\varphi$ on $\rr^n$ in the Alexandrov a.e. sense and ${\mathcal A}(x,\mu_t)=\nabla V(x)+\nabla W*\mu_t(x)$. Moreover, since $\nabla W$ is odd, the term involving $W$~is
\begin{eqnarray*}
&& \int_{\rr^n} (\nabla W*\nu_t(\nabla\phi_t(x)) -\nabla W*\mu_t(x)) \cdot (\nabla\phi_t(x)-x) \, d\mu_t(x)\\
&& \!\!\!\!\!\!\!\!\!=\iint_{\rr^{2n}} (\nabla W(\nabla\phi_t(x)-\nabla\phi_t(y))-\nabla W(x-y)) \cdot (\nabla\phi_t(x)-x) \, d\mu_t(x)d\mu_t(y)\\
&&\!\!\!\!\!\!\!\!\! =\frac12\iint_{\rr^{2n}} (\nabla W(\nabla\phi_t(x)-\nabla\phi_t(y))-\nabla W(x-y)) \cdot (\nabla\phi_t(x)-\nabla\phi_t(y)-(x-y)) \, d\mu_t(x)d\mu_t(y).
\end{eqnarray*}

We summarize as follows:

\begin{eprop}[\cite{ambrosio-gigli-savare}]
\label{derWgran}
If $(\mu_t)_t$ and $(\nu_t)_t$ are two solutions to~\eqref{eq-gran}, then for a.e.  $t >0$, 
$$
\frac12 \frac{d}{dt}W_2^2(\mu_t,\nu_t) \leq -J_{V,W}(\nu_t \vert \mu_t)
$$
where, for $\nu=\nabla\phi\#\mu$ (and $\Delta\phi$  the trace of the Hessian of $\varphi$ in the Alexandrov sense), 
\begin{multline}
 \label{eq-jgran}
 J_{V,W}(\nu|\mu)
\!=\!\!\! \int_{\rr^n} \Big( \Delta\phi(x) +\Delta\phi^*(\nabla\phi(x))-2n + (\nabla V(\nabla\phi(x))-\nabla V(x)) \cdot (\nabla\phi(x)-x) \Big) \,d\mu(x) \\
+\frac12\iint_{\rr^{2n}}  (\nabla W(\nabla\phi(x)-\nabla\phi(y))-\nabla W(x-y)) \cdot (\nabla\phi(x)-\nabla\phi(y)-(x-y))d\mu(x)d\mu(y).  
\end{multline}
\end{eprop}

For $t>0$ we can expect the solutions to have smooth densities, and to have equality in Proposition~\ref{derWgran}, but we shall be content with the inequality (see~\cite{calvez-carrillo10}).

\medskip

Considering the dissipation of the distance between two solutions provides simple alternative proofs of contraction properties in Wasserstein distance derived in~\cite{cmcv-06,cgm-08}. 
For that purpose we first notice that given $\mu$ and $\nu$ absolutely continuous with respect to the Lebesgue measure, and $\nabla\phi\#\mu=\nu$, then  $\Delta\phi+\Delta\phi^*(\nabla\phi)-2n\geq0$ $\mu$ a.e. (see for example ~\cite[Th. 1.5]{lisini} and~\cite[Lem.~2.5]{BGG11}). 
This inequality says that the diffusion part of the equation always contracts two solutions, as it is classical for the pure heat equation. Then:
\begin{itemize}
\item Suppose that  $V$ and $W$ are respectively $\alpha$ and $\beta$-convex with $\alpha \in \rr$ and $\beta \leq 0$.
Then 
the term involving $V$ in~\eqref{eq-jgran} is bounded from below by $\alpha \, W_2^2(\mu, \nu)$ and the term involving~$W$~by
$$
\!\!\!\!\!\!\!\!\!\frac\beta2\iint  \big\vert \nabla\phi(x)-\nabla\phi(y)-(x-y) \big\vert^2d\mu(x)d\mu(y)
\\
 = \beta W_2^2(\mu, \nu) - \beta \Big\vert \int (\nabla \phi - x) \, d\mu \Big\vert^2
 \geq \beta W_2^2(\mu, \nu)
$$
since $\beta \leq 0$. Hence, for two solutions $(\mu_t)_t$ and $(\nu_t)_t$ of~\eqref{eq-gran}  and almost all $t\geq0$ 
 $$
 \frac12 \frac{d}{dt}W_2^2(\mu_t,\nu_t)\leq -(\alpha + \beta) \, W_2^2(\mu_t, \nu_t).
 $$
Then by the Gronwall lemma we recover the contraction property of~\cite[Th. 5]{cmcv-06}:
 \begin{equation}\label{contr2exp}
 W_2(\mu_t,\nu_t)\leq e^{-(\alpha + \beta)t} \, W_2(\mu_0, \nu_0), \qquad t \geq 0.
 \end{equation}

\item Suppose that $W$ is convex and that there exist $p,C>0$ such that for all $\varepsilon>0$
 \begin{equation}\label{Vconvdege}
(\nabla V(y)- \nabla V(x)) \cdot (y-x)\ge C\,\varepsilon^p(|y-x|^2-\varepsilon^2), \qquad x, y \in \rr^n.
\end{equation}
Then, by the same argument,
$$
 \frac12 \frac{d}{dt}W_2^2(\mu_t,\nu_t) \leq - \frac{C}{2} \varepsilon^p(2W_2^2(\mu_t,\nu_t)-\varepsilon^2).
$$
We optimize in $\varepsilon$ and integrate to recover the polynomial contraction of~\cite[Th. 6]{cmcv-06} 
 \begin{equation}\label{contr2pol}
W_2(\mu_t,\nu_t)\le \left(W_2^{-p}(\mu_0,\nu_0)+c t\right)^{-1/p}, \qquad t \geq 0.
\end{equation}

\item Suppose that $V$ and $W$ are respectively $\alpha$ and $\beta$-convex with $\alpha \in \rr$ and $\beta \geq 0$. Then, again by the same argument, the contraction result~\eqref{contr2exp} holds for any two solutions with same center of mass, that is, such that in $\rr^n$
$
\int_{\rr^n} xd\mu_t =\int_{\rr^n} x d\nu_t 
$
 for all $t\geq0$. This was also proved in~\cite[Th. 5]{cmcv-06}. 

\item Suppose that $V$ is convex and that~\eqref{Vconvdege} holds for $W$ instead of $V$.
If moreover the center of mass of each solution is conserved, that is, if $\int_{\rr^n} x  d\mu_t=\int_{\rr^n}  x  d\mu_0$ for all $t \geq 0$ (this is the case for instance if $V=0$) then the polynomial contraction~\eqref{contr2pol} holds for any two solutions with same (initial) center of mass, recovering~\cite[Th. 6]{cmcv-06} and~\cite[Th.~4.1]{cgm-08}. 
\end{itemize}

In the first case with $\alpha + \beta >0$  the bound~\eqref{contr2exp} ensures the existence of a unique stationary solution to~\eqref{eq-gran} in $P_2(\rr^n)$, and the exponential  convergence of all solutions to it. In the third case with $\alpha + \beta >0$, and if moreover the center of mass is preserved by the evolution, then for any $m \in \rr^n$ this ensures the existence of a unique stationary solution to~\eqref{eq-gran} in $P_2(\rr^n)$  with center of mass $m$, and the exponential  convergence to it of all solutions with (initial) center of mass $m$.

\medskip

 The following two sections are devoted to the obtention of explicit exponential rates of convergence of solutions to~\eqref{eq-gran} in non uniformly convex or even non convex cases, having in mind the degenerately convex potentials of~\cite{BCCP} and the double well potentials of~\cite{tugaut}.

\section{Influence of the interaction potential}\label{sec2}

In this section we study the case when $W$ brings the convergence.

\subsection{No exterior potential}

We first assume that $V=0$.  Then the evolution preserves the center of mass, and a solution $\mu_t$ should converge to a stationary solution $\mu_{\infty}$ only if the initial datum $\mu_0$ and $\mu_{\infty}$ have same center of mass, since
$$
\int_{\rr^n} x \, d\mu_{\infty} (x) - \int_{\rr^n} x \, d\mu_0(x) = \int_{\rr^n} x \, d\mu_{\infty}(x) - \int_{\rr^n} x \, d\mu_t(x)
$$
should converge to $0$: for instance it is bounded by $W_2(\mu_t, \mu_{\infty})$. We could also assume that $V \neq 0$, but that the center of mass is fixed by the evolution, which is all we use. But to simplify the statements we assume $V=0$.

When $W$ is degenerately convex, with a pointwise degeneracy, for instance $W(x) = \vert x \vert^{2+ \varepsilon}$ with $\varepsilon >0$, then the contraction property holds only with polynomial decay rate, see the last example in Section~\ref{sec1}. Then in~\cite{cmcv-03} the authors proved an exponential convergence to equilibrium, but not with a uniform decay rate, but rather depending on the free energy $F$ of the initial datum. In this section we prove a {\it uniform} exponential convergence for such potentials.

\begin{ethm}\label{thm:main}
Let $V=0$ and $W$ a $\mathcal C^2$ convex map on $\rr^n$ for which there exist $R$ and $K>0$ such that
$$
\nabla^2 W(x) \geq K \quad  \textrm{if} \quad \vert x \vert \geq R.
$$
Then for all $m \in\rr^n$ there exists a unique stationary solution $\mu_{\infty}^m \in P_2(\rr^n)$ to~\eqref{eq-gran} with center of mass $m$; moreover there exists a positive  constant $C$ such that all solutions $(\mu_t)_t$ to~\eqref{eq-gran}, for an initial datum with center of mass $m$, converge to $\mu_{\infty}^m$ according to
$$
W_2 (\mu_t, \mu_{\infty}^m) \leq e^{-Ct} W_2(\mu_0, \mu_{\infty}^m), \qquad t \geq 0.
$$
\end{ethm}

\begin{eproof}
The existence of a stationary solution $\mu_{\infty}^0 \in P_2(\rr^n)$ with center of mass $0$ and a positive density satisfying $\mu_{\infty}^0(x) = Z^{-1} e^{-W * \mu_{\infty}^0(x)}
$
is given by Proposition~\ref{prop:H}, proof of $i.$, with any $b < K/2$; here $Z$ is the normalizing constant.  Then $\mu_{\infty}^m = \mu_{\infty}^0(\cdot -m)$ is a stationary solution with center of mass $m$. Now Proposition~\ref{prop:granW} and Remark~\ref{rem:hypU}, $ii.$ below ensure the convergence estimate to $\mu_{\infty}^m$ since $\mu_{\infty}^m = e^{-U} /Z$ with $U = W * \mu_{\infty}^m$ convex and bounded from below. Uniqueness follows. 
\end{eproof}

\begin{eprop}\label{prop:granW}
Let $W$ be a $\mathcal C^2$ convex map on $\rr^n$ for which there exist $R$ and $K>0$ such that
$$
\nabla^2 W(x) \geq K \quad  \textrm{if} \quad \vert x \vert \geq R.
$$
Let $\mu \in P_2(\rr^n)$ have a continuous density $e^{-U}$ for which there exists $M$ such that 
\begin{equation}\label{hypU}
\sup_{\vert x-y \vert \leq 2R}\,\,\sup_{z \in [x,y]} \{U(z) - U(x)- U(y) \}\leq M.
\end{equation}

Then there exists an explicit positive constant $C$, depending only on $K, R$ and $M$, such that 
$$
C  \, W_2^2(\nu, \mu) \leq J_{0,W}(\nu \vert \mu)
$$
for all measures $\nu$ with same center of mass as $\mu$.
\end{eprop}

\begin{erem}\label{rem:hypU}
Hypothesis~\eqref{hypU} on  $U$ holds on any of the following two instances:

\begin{enumerate}
\item $U$ is $\mathcal C^1$ and $U(x) - \displaystyle 2R \sup_{\vert x-y \vert \leq R} \vert \nabla U(y) \vert \geq -M$ for all $x \in \rr^n$

\item $U$ is $\mathcal C^2$, $\nabla^2 U(x) \geq \alpha (x)$ with $\alpha(x) \leq 0$ and $U(x) + 2R^2 \displaystyle \inf_{\vert x-y \vert \leq R} \alpha (y)
\geq -M $ for all $x$; for example, $U$ is $\mathcal C^2$ and bounded from below and $\nabla^2 U(x) \geq \alpha$ for all $x$ and a constant $\alpha$.

\end{enumerate}

For $ii.$ for instance, assume that the sup of $U$ on $[x,y]$ is achieved at $z= tx + (1-t) y$ with $t \in ]0,1[$. Then $\nabla U(z) \cdot (y-x) =0$, so that
$$
U(x) - U(z) = \int_0^1 (1-s) \nabla^2 U(z+ s(x-z)) (x-z) \cdot (x-z) \, ds 
\geq \inf_{\vert Y - \frac{x+y}{2} \vert \leq R} \alpha (Y) \, \frac{(1-t)^2}{2} \vert x-y\vert^2
$$
and similarly
$$
U(y) - U \big( \frac{x+y}{2} \big) \geq U(y) - U(z) 
\geq \inf_{\vert Y - \frac{x+y}{2} \vert \leq R} \alpha (Y) \, \frac{t^2}{2} \vert x-y\vert^2.
$$
Hence, for $\vert x-y \vert \leq 2R$,
$$
U(z) - U(x) - U(y) \leq -U\big( \frac{x+y}{2} \big) - 2 R^2 \inf_{\vert Y - \frac{x+y}{2} \vert \leq R} \alpha (Y) \leq M.
$$

\end{erem}


\medskip

\noindent {\emph{\textbf{Proof of Proposition~\ref{prop:granW}}}\\\proofbegin~
Let  $\phi$ be a strictly convex function on $\rr^n$ (with $\nu = \nabla \phi \# \mu$) such that $ \int_{\rr^n} \nabla \phi \, d\mu = \int_{\rr^n} x \, d\mu$.

First observe that
$$
\iint_{\rr^{2n}} \vert \nabla \phi(x) - \nabla \phi(y) - (x-y) \vert^2 \, d\mu(x) \, d\mu(y)
= 2 \int_{\rr^n} \vert \nabla \varphi(x) - x \vert^2 \, d\mu(x) 
$$
since, by assumption on $\phi$,  the difference is
$$
2 \; \Big\vert \int_{\rr^n} (\nabla \phi(x) - x) \, d\mu(x) \Big\vert^2 = 0.
$$
Then, by~\cite[Lem.~5.1]{BGG11},
\begin{equation}
\label{eq-lemma}
(\nabla W (x) - \nabla W (y))  \cdot (x-y) \geq \frac{K}{3}\vert x-y \vert^2
\end{equation}
if $\vert x \vert \geq 2R$ or  $\vert y \vert \geq 2R$. In view of this result we let
$$
X = \{(x,y) \in \rr^{2n}; \vert x-y \vert \leq 2R, \vert \nabla \phi(x) - \nabla \phi(y) \vert \leq 2R \}.
$$

{\bf 1.} First of all, by convexity of $W$ and~\eqref{eq-lemma},
  \begin{eqnarray*}
&& \int_{\rr^{2n}} (\nabla W (\nabla \phi(x) - \nabla \phi (y))- \nabla W(x-y)) \cdot  (\nabla \phi(x) - \nabla \phi(y) -(x-y)) \,d\mu(x) d\mu(y) 
 \\
&&\geq   \int_{\rr^{2n} \setminus X} (\nabla W (\nabla \phi(x) - \nabla \phi(y))- \nabla W (x-y)) \cdot (\nabla \phi(x) - \nabla \phi(y) -(x-y)) \,d\mu(x) d\mu(y) 
\\
&&\geq 
\frac{K}{3}  \int_{\rr^{2n}\setminus X}\vert \nabla \phi(x) - \nabla \phi (y) - (x-y) \vert ^2 \, d\mu(x) d\mu(y).
 \end{eqnarray*}

{\bf 2.} Then for all $x$ and $y$, written as $y = x + r \, \theta$ with $r \geq 0$ and $\theta \in \mathbb S^{n-1}$, 
$$
\nabla \phi (y) -  \nabla \phi (x) - (y-x) = \int_0^1 [\nabla^2 \phi(x + r \,t\, \theta) - I] \, r \, \theta \, dt.
$$
We let $H = \nabla^2 \phi(x + r \, t\, \theta)$ and write $H - I = [H^{1/2} - H^{-1/2}] H^{1/2}$, so that
$$
\vert [H - I] \theta \vert \leq \Vert H^{1/2} - H^{-1/2} \Vert \; \;  \vert H^{1/2} \theta \vert . 
$$
Hence
$$
\vert  \nabla \phi (y) -  \nabla  \phi (x) - (y-x) \vert^2
\leq
 r \int_0^1  \Vert H^{1/2} - H^{-1/2} \Vert^2 \, e^{-U(x+r \, t\, \theta)} dt   \int_0^1 \vert H^{1/2} \theta \vert^2 e^{U(x+r \, t\, \theta)} \, r \, dt
$$
 by the Cauchy-Schwarz inequality.
 On the one hand, letting $D = \Delta\phi  +\Delta\phi^*(\nabla\phi )-2n$,
$$
\Vert H^{1/2} - H^{-1/2} \Vert^2 
=\sup_x \frac{([H-2I+H^{-1}]x) \cdot x}{\vert x \vert^2}
\leq
  trace (H-2I+H^{-1}) = D(x + r \, t\, \theta).
$$
On the other hand
$$
 \int_0^1 \vert H^{1/2} \theta \vert^2 \, r \, dt
=
 \int_0^1 \nabla^2 \phi(x+r \, t\,\theta) (r\theta) \cdot \theta \, dt \\
 =
 (\nabla \phi(y) - \nabla \phi(x)) \cdot \theta  \leq 2 R 
$$
 if $(x,y) \in X$.
Hence 
$$
\vert  \nabla \phi (y) -  \nabla  \phi (x) - (y-x) \vert^2
\leq
4 R^2 \sup_{z \in [x,y]} e^{U(z)} \int_0^1 D(x + r\, t \, \theta) e^{-U(x+r \, t\, \theta)} dt
$$
for all $(x,y) \in X$, so that
\begin{multline*}
\iint_X \vert  \nabla \phi (y) -  \nabla  \phi (x) - (y-x) \vert^2 \, d\mu(x) d\mu(y)
\\
\leq
4 R^2 \int_{\rr^n} e^{-U(x)} \, dx  \int_{\vert y-x  \vert \leq 2R} dy \, \sup_{z \in [x,y]} e^{U(z)} e^{-U(y)} \int_0^1 D(x + t (y-x)) \, e^{-U(x+t (y-x))} dt
\\
\leq 4 R^2 e^M \int_{\rr^n} dx \int_{\vert y-x  \vert \leq 2R} dy \int_0^1 D(x + t (y-x)) \, e^{-U(x+t (y-x))} dt
\end{multline*}
by~\eqref{hypU}.  Now, for fixed $t \in [0,1]$, the change of variables $(x,y) \mapsto (v,u) = (x+t(y-x), y-x)$ has unit Jacobian, so this is equal to
$$
4 R^2 e^M \int_0^1 dt \int_{\rr^n} dv \int_{\vert u \vert \leq 2R} du  \, D(v) e^{-U(v)} 
 = c \int_{\rr^n}  D(v) \, d\mu(v)
 $$
for a constant $c = c(R,M,n) = 4 R^{2+n} e^M c_n$ where $c_n$ is the volume of the unit ball in $\rr^n$.

 \smallskip
 
 {\bf 3.} Collecting the terms in {\bf 1.} and {\bf 2.} concludes the proof with $C = 2 (\frac{3}{K} + 4 c_n R^{2+n} e^M )^{-1}.$~\proofend

\subsection{In presence of an exterior potential}

 We saw in Section~\ref{sec1} how an exterior potential $V$ can induce the convergence of all solutions to a unique equilibrium, and not only to the unique equilibrium with same center of mass as the initial datum of the solution to be considered.
 
If  W strictly convex (but at $0$), and uniformly at infinity, and if $V$ is strictly convex (but at $0$), then polynomial convergence holds to a unique equilibrium $\mu_{\infty}$ (see Section~\ref{sec1} and~\cite[Th.~2.3]{cmcv-03}), and even exponential convergence, but with a rate depending on the free energy $F$ of the initial datum (see~\cite[Th.~2.5]{cmcv-03}). Following Theorem~\ref{thm:main} (for $V=0$), one may wonder whether this convergence is actually uniform in the initial datum, given by
 $$
 W_2(\mu_t, \mu_{\infty}) \leq e^{-Ct} W_2(\mu_0, \mu_{\infty}), \qquad t \geq 0
 $$
 for all solution $(\mu_t)_t$. But, according to Section~\ref{sec1}, this estimate is based on the inequality
 \begin{equation}\label{WJ}
 C \, W_2^2(\nu, \mu) \leq J_{V,W} (\nu| \mu)
 \end{equation}
 for the measure $\mu = \mu_{\infty}$ and all measures $\nu$; and this inequality does not hold if $V$ is only assumed to be strictly convex. {\it For instance:}
 
 \begin{elem}\label{lem:nonWj}
Let $\mu \in P_2(\rr)$ and $V$ be a $\mathcal C^2$ map on $\rr$, with $V''$ bounded and $\displaystyle V'' \to_{+ \infty} 0$. Then there is no constant $C>0$ such that~\eqref{WJ} hold for all $\nu$.
 \end{elem}

  \proofbegin~We prove that~\eqref{WJ} does not hold for the translations $\nu = \phi' \# \mu$ where $\phi'(x) = x + M$ where $M \to + \infty$, that is, that there is no $C>0$ such that
\begin{equation}\label{eq:nonWJ}
  CM \leq  \int_{\mathbb R}  (V'(x+M) - V'(x)) \, d\mu(x)
 \end{equation}
 for all $M>0$.
  For that, we let $R$ to be fixed later on, and bound the right-hand side in~\eqref{eq:nonWJ} by
 $$
 \int_{\rr} \vert V'(x) \vert \, d\mu(x) +  \int_{\infty}^{-R} \vert V'(x+M) \vert \, d\mu(x) +  \int_{-R}^{+ \infty} \vert V'(x+M) \vert \, d\mu(x).
 $$

 First of all, since $\vert V'' \vert \leq A$, then $\vert V'(x) \vert \leq \vert V'(0) \vert + A \vert x \vert$ so the first integral is finite (uniformly in $M$), and the second one is bounded by
 $$
 \alpha M \int_{\infty}^{-R} d\mu(x) + \int_{\infty}^{-R} ( \vert V'(0) \vert + A \vert x \vert) \, d\mu(x).
 $$
 Now, for fixed $\varepsilon >0$, we take $R$ such that this is bounded by $(M+1) \varepsilon$ for all $M$. Then we take $M_0$ such that $\vert V''(x) \vert \leq \varepsilon$ for $x \geq M_0$. For all $M \geq M_0 + R$ the third integral is bounded by
 $$
 \int_{-R}^{+ \infty} \big( \vert V'(M_0) \vert + \varepsilon (x+M-M_0) \big)\, d\mu(x) 
 \leq 
\vert V'(M_0) \vert + \varepsilon \Big( M + M_0 + \int_{\rr} \vert x \vert \,  d\mu(x)  \Big).
$$
Collecting all terms we conclude that the full right-hand side in~\eqref{eq:nonWJ} is $\leq 4 \varepsilon$ for large $M$.
  \proofend
  
\medskip

Lemma~\ref{lem:nonWj} only gives an instance of condition on $V$ for~\eqref{WJ} not to hold. For example, the assumption $V''$ bounded can be replaced by 
the doubling condition~\eqref{eq-doubling}  for $V'$  and $ \int \vert V' \vert \, d\mu <~\infty$. Similarly, the assumption $V'' \to_{+ \infty} 0$ can be replaced by $\int_0^{+ \infty} \vert V'' \vert \, d\mu <~\infty$: in this case we~use
$$
\Big\vert 
\int_{-R}^{+ \infty} 
\! 
(V'(x+M) - V'(x)) \, d\mu(x)
\Big\vert
=
\Big\vert 
\!
\int_{-R}^{+ \infty} 
\! \! \!
\int_{x}^{x+M} 
\! \! \!
V''(t) \, dt  \, d\mu(x)
\Big\vert
  \leq 
\! \! 
\int_{-R}^{+ \infty} 
\! \! \!
\int_{\rr} \vert V''(t) \vert \, dt  \, d\mu
\leq 
C.
$$

  \bigskip
  
  Hence we can not expect a uniform rate of convergence to equilibrium for degenerately convex potentials. Our method is however able to recover an exponential convergence with a rate depending on the initial datum, as in~\cite[Th. 2.5]{cmcv-03}:
\begin{ethm}
Assume that $V$ is convex on $\rr^n$ with $\nabla ^2 V(x)$ definite positive on $\vert x \vert \geq R$ and $\int e^{-V}dx<\infty$, and that $W$ is convex with $\nabla ^2 W(x) \geq K$ for $\vert x \vert \geq R$. Then there exists a unique stationary solution $\mu_{\infty} \in P_2(\rr^n)$ to~\eqref{eq-gran}. Moreover for all $M$ there exists a positive constant~$C$ such that for all solutions $(\mu_t)_t$ with $ \int_{\rr^n} \vert x \vert^2 d\mu_0(x) \leq M$
$$
W_2(\mu_t, \mu_{\infty}) \leq e^{-Ct} \, W_2(\mu_0, \mu_{\infty}), \quad t \geq 0.
$$
\end{ethm}

\proofbegin~ 
First, Proposition~\ref{prop:H}, $ii.$ ensures the existence of a stationary measure $\mu_{\infty} \in P_2(\rr^n)$ which has a density satisfying $
\mu_{\infty}(x) = Z^{-1} e^{-V(x) - W * \mu_{\infty}(x)};$
here $Z$ is the normalizing constant. We just mention that the assumptions on $V$ are satisfied by~\cite[Lem. 2.2]{bbcg08} for instance.

\medskip

Then, by direct estimates on the propagation of the second moment, for all solutions $(\mu_t)_t$ with $\int_{\rr^n} \vert x \vert^2 d\mu_0(x) \leq M$ there is a constant $N$, depending only on $V, W$ and $M$  such that
$$
\sup_{t \geq 0}  \int_{\rr^n} \vert x \vert^2 d\mu_t(x) \leq N.
$$
 
 \medskip
Moreover, for $\nu = \nabla \varphi \# \mu_{\infty}$ with $\int_{\rr^n} \vert x \vert^2 d\nu(x) \leq N$, we first write
$$
\int_{\rr^n}\! \! \! \! \vert \nabla \varphi(x) - x \vert^2 \, d\mu_{\infty} \! = \!  \Big\vert \! \! \int_{\rr^n} \!\! \! \! (\nabla \phi(x) - x) \, d\mu_{\infty}\Big\vert^2
\! + \frac{1}{2}   \! \iint_{\rr^{2n}} \! \! \! \!\! \vert \nabla \phi(x) - \nabla \phi(y) - (x-y) \vert^2 \, d\mu_{\infty}(x)  d\mu_{\infty}(y).
$$
By Proposition~\ref{prop:vstrictcvx} below, applied with the constant $N$ and the measure $\mu_{\infty}$, 
$$
 \Big\vert \! \! \int_{\rr^n} \!\! \! \! (\nabla \phi(x) - x) \, d\mu_{\infty}\Big\vert^2 \leq \frac{1}{2} \int_{\rr^n} \vert \nabla \varphi(x) - x \vert^2 \, d\mu_{\infty}(x) + C \, J_{V, 0} (\nu \vert \mu_{\infty}).
$$
Then, by the proof of Proposition~\ref{prop:granW}, there exists $C_1$, depending only on $V$ and $W$, such that
$$
  \iint_{\rr^{2n}} \vert \nabla \phi(x) - \nabla \phi(y) - (x-y) \vert^2 \, d\mu_{\infty}(x) \, d\mu_{\infty}(y)
\leq C_1 J_{0,W}(\nu \vert \mu_{\infty}).
$$
Hence there exists a new positive constant $C$, depending only on $V, W$ and $M$, such that
$$
C W_2^2(\mu_t, \mu_{\infty}) \leq J_{V,W}(\mu_t \vert \mu_{\infty})
$$
for all $t$. This proves the estimate on the convergence to $\mu_{\infty}$ again by Proposition~\ref{derWgran} and the Gronwall lemma. Uniqueness of the stationary solution in $P_2(\rr^n)$ follows.
  \proofend

\begin{eprop}\label{prop:vstrictcvx}
Let $V$ be a $\mathcal C^2$ convex map on $\rr^n$ with $\nabla ^2 V(x)$ definite positive on $\vert x \vert \geq R$, and $d\mu(x) = e^{-U(x)} \, dx$ be a probability measure on $\rr^n$ with $U$ continuous. Then for all $N$ there exists a constant $C$ such that for all $\mathcal C^2$ strictly convex map $\varphi$ on $\rr^n$ with $ \int_{\rr^n} \! \! \vert \nabla \varphi(x) \vert^2 \, d\mu \leq N$
$$
  \Big\vert \int_{\rr^n} \nabla \phi(x) \, d\mu(x)  - \int x \, d\mu(x) \Big\vert^2
  \leq
\frac{1}{2}  \int_{\rr^n} \vert \nabla \phi(x) -x \vert ^2 d\mu(x) 
 +
  C \, J_{V,0} (\nabla \varphi \# \mu \vert \mu). 
$$
\end{eprop}

\proofbegin~
Let $S  \geq 3R$ to be fixed later on. Since $V$ is $\mathcal C^2$ and $\nabla ^2 V(x)$ is definite positive on the compact set $R \leq \vert x \vert \leq S$, there exists $K = K(S)>0$ such that $\nabla^2 V (x) \geq K$ for all $R \leq \vert x \vert \leq S$. Then, following~\cite[Lem.~5.1]{BGG11},
$$
(\nabla V(y) - \nabla V(x)) \cdot (y-x) \geq \frac{K}{3} \vert x-y \vert^2
$$
if $\vert x \vert \leq S$, $\vert y \vert \leq S$ and if $\vert x \vert \geq 2R$ or $\vert y \vert \geq 2R$; indeed one only need to take into account the values of $\nabla^2 V$ on the ball of radius $S$. 

Then we let $\varphi$ be a given $\mathcal C^2$ strictly convex map on $\rr^n$ and let
$$
X_1 = \{x \in \rr^n, \vert x \vert \leq S, \vert \nabla \varphi(x) \vert \leq S, \vert x \vert \geq 2R \; \textrm{or} \; \vert \varphi(x) \vert \geq 2R \}.
$$

\medskip

{\bf 1.} First of all, by convexity of $V$, the above remark and the Cauchy-Schwarz inequality,
\begin{eqnarray*}
\int_{\rr^n} (\nabla V(\nabla \varphi) - \nabla V) \cdot (\nabla \varphi-x) d\mu
&\geq& 
\int_{X_1} (\nabla V(\nabla \varphi) - \nabla V) \cdot (\nabla \varphi-x) d\mu
\\
&\geq& 
\frac{K}{3} \int_{X_1} \vert \nabla \varphi - x \vert^2 d\mu
\geq 
\frac{K}{3} \Big \vert \int_{X_1} (\nabla \varphi - x)  d\mu \Big\vert^2.
 \end{eqnarray*}

\medskip

{\bf 2.} Then, on $\rr^n \setminus X_1$, and letting $X_2 = \{x \in \rr^n, \vert x \vert \leq 2R, \vert \nabla \varphi(x) \vert \leq 2R \}$,
$$
\Big \vert \int_{\rr^n \setminus X_1} \! \! (\nabla \varphi - x)  d\mu \Big\vert^2
\leq
3 \Big( \int_{\vert x \vert \geq S} \vert \nabla \varphi -x \vert d\mu \Big)^2
+
3 \Big( \int_{\vert \nabla \varphi(x) \vert \geq S} \vert \nabla \varphi -x\vert d\mu \Big)^2
+
3 \Big( \int_{X_2} \! \! \vert \nabla \varphi -x \vert d\mu \Big)^2
 $$

By the Cauchy-Schwarz and Markov inequalities, the first term is bounded from above by
$$
\int_{\vert x \vert \geq S} \vert \nabla \varphi -x \vert^2 d\mu \; \; \; \mu [x , \vert x \vert \geq S] \leq \int_{\rr^n} \vert \nabla \varphi -x \vert^2 d\mu \; \frac{1}{S^2} \int_{\rr^n} \vert x \vert^2 d\mu(x)
$$
and the second one by
$$
\int_{\rr^n} \vert \nabla \varphi -x \vert^2 d\mu \; \frac{1}{S^2} \int_{\rr^n} \vert \nabla \varphi (x) \vert^2 d\mu(x).
$$

Then, following the proof of~\cite[Prop. 3.5]{BGG11}, there exists a constant $C$, depending on $V$ and $\mu$ only on the ball of radius $3R$, such that
$$
 \int_{X_2} \vert \nabla \phi -x \vert ^2 d\mu
 \leq
  C \int_{\vert x \vert \leq 3R} \big( \Delta\phi+\Delta\phi^*(\nabla\phi)-2n  + (\nabla V(\nabla\phi)-\nabla V) \cdot (\nabla\phi-x) \big)\,d\mu.
$$
This is in turn bounded by the corresponding integral on $\rr^n$, which is $J_{V,0} (\nabla \phi \# \mu \vert \mu)$, since both terms in the integrand are nonnegative.

\medskip

{\bf 3.} Collecting all terms we obtain
 $$
 \Big\vert \int_{\rr^n} \! \! \nabla \phi \, d\mu - \int_{\rr^n} \! \! x \, d\mu \Big\vert^2
  \leq
\frac{6}{S^2}  \int_{\rr^n} \! \! \vert \nabla \phi -x \vert ^2 d\mu \Big[ \int_{\rr^n} \! \! \vert x \vert^2 d\mu +  \int_{\rr^n} \! \! \vert \nabla \varphi \vert^2 d\mu \Big]
 +
\Big(\frac{6}{K} + 6 \, C \Big) J_{V,0} (\nabla \phi \# \mu \vert \mu). 
$$

Then we let $S = \max \Big\{ 3R, \sqrt{12} \Big[ \displaystyle \int \vert x \vert^2 d\mu + N \Big] \Big\}$ so that
$$
\frac{6}{S^2} \Big[ \int \vert x \vert^2 d\mu +  \int \vert \nabla \varphi \vert^2 d\mu \Big] \leq \frac{1}{2}
$$
if $\int \vert \nabla \varphi  \vert^2 d\mu \leq N$, concluding the proof with a $C$ depending on $V, \mu$ and $M$ through $K(S)$.
\proofend
\section{Non convex examples}\label{sec3}

In this section we deal with potentials $V$ and $W$ for which the convergence rate to equilibrium is driven by $V$ rather than by $W$. 
Our first result is more  qualitative rather than quantitative.

\begin{ethm}
\label{thm-last}
Assume that $V$ and $W$ are $\mathcal C^2$ convex maps and that there exist $R\ge 0$ and $K>0$ such that for all $|x|\ge R$,
$$
\nabla^2 V(x) \ge K.
$$
 Then  there exists a unique stationary solution $\mu_{\infty} \in P_2(\rr^n)$ to~\eqref{eq-gran}, and a constant $C$ such that for all solution $(\mu_t)_t$ of~\eqref{eq-gran},
$$
W_2 (\mu_t, \mu_{\infty})\leq e^{-Ct}W_2 (\mu_0, \mu_{\infty}), \qquad t \geq 0.
$$
\end{ethm}

\bigskip

In the first section (second example) we saw that only polynomial decay in contraction is known in this context, and only when the convexity degenerates at some points, for instance for $V(x)=|x|^4$.

\bigskip

\begin{eproof}~
Existence of a stationary solution $\mu_{\infty}$ in $P_2(\rr^n)$ which has a positive density satisfying
$\mu_{\infty} = Z^{-1} e^{-V - W * \mu_{\infty}}$
is given by Proposition~\ref{prop:H}, $iii$, with any $a<K$ and $-a < b < 0$.  Then, by
 \cite[Prop.~3.5]{BGG11},  there exists $C>0$ such that 
$$
C\,W_2^2(\mu_t,\mu_\infty)\le  J_{V,0}(\mu_t\vert   \mu_{\infty})
$$
for all solution $(\mu_t)_t$. Moreover $ W$ is convex, so $J_{V,0} \leq J_{V,W}$. This proves the convergence bound by Proposition~\ref{derWgran}. Uniqueness of the stationary solution in $P_2(\rr^n)$ follows.
\end{eproof}

\medskip

\begin{erem}
The case of a double well potential for $V$ is considered by J.~Tugaut in~\cite{tugaut2012, tugaut}, where the long time behavior is studied by a compactness argument, hence without rate.   Let us now explain  how Theorem~\ref{thm-last} extends to this case, for instance for $V^\varepsilon(x)=x^4-\varepsilon x^2$ and $W(x)=|x|^3$ in $\rr$. 

 First of all, a stationary solution, solution of $\mu_\infty^\varepsilon=e^{-V^\varepsilon-\mu_{\infty}^\varepsilon*W}/Z^\varepsilon$, exists by Proposition~\ref{prop:H}, $iii$. Then one can then easily build a cut-off function $\psi$ such that $V^\varepsilon\psi$ is $\mathcal C^2$, convex, satisfies $(V^\varepsilon\psi)'' \geq K>0$ outside a centered ball, uniformly in $\varepsilon \in [0,1]$,  and is such that $\|(V^\varepsilon(1-\psi))''\|_\infty$ converges to 0 as $\varepsilon\to0$.
Then, by  \cite[Prop.~3.5]{BGG11}, the measure $\mu_\infty^\varepsilon$ satisfies a $WJ_{V^\varepsilon \psi,0}$ inequality with a constant $C >0$ uniformly in $\varepsilon\in [0,1]$ (here we use that $\int Wd\mu_\infty^\varepsilon$ and $Z^\varepsilon$ are bounded uniformly in $\varepsilon$). Now the perturbation proposition \cite[Prop.~3.8]{BGG11} ensures that $\mu_\infty^\varepsilon$ satisfies a $WJ_{V^\varepsilon,0}$ inequality, for sufficiently small $\varepsilon$, hence a $WJ_{V^\varepsilon,W}$ inequality since  $W$ is convex. Here we say that a measure $\mu$ satisfies a $WJ_{V,W}$ inequality is the inequality~\eqref{WJ} holds for a positive constant $C$ and all $\nu$. 

The smallness condition on $\varepsilon$ is necessary since, according to~\cite{tugaut2012, tugaut}, there exist several stationary solutions for large $\varepsilon$. 
\end{erem}

\medskip

 The following theorem provides the first examples of exponential convergence to equilibrium  for the granular media equation, with both potentials non convex.

\smallskip

\begin{ethm}
Assume that
\begin{itemize}
\item $e^{-V} \in P_2(\rr^n)$ and there exist $\alpha \in\R$ and $C>0$ such that $\nabla^2 V\ge \alpha $, and  for all $\nu$
\begin{equation}
\label{eq-defn}
 W_2^2(\nu, e^{-V})\leq{\frac1C \,J_{V,0}(\nu| e^{-V})};
\end{equation}
 
\item there exist $K \geq 0$ and $\beta \leq 0$ such that $\sup |W| \le K$ and $\nabla^2 W\ge\beta $.
\end{itemize}
Then there exists a unique stationary solution $\mu_{\infty} \in P_2(\rr^n)$ to~\eqref{eq-gran}. Moreover, for all solution $(\mu_t)_t$ of ~\eqref{eq-gran}, and with $\tilde C= (C- \alpha) e^{-2K} + \alpha + \beta,$
$$
W_2 (\mu_t, \mu_{\infty})\leq e^{-\tilde Ct}W_2 (\mu_0, \mu_{\infty}), \quad t\geq 0.
$$
\end{ethm}

Assumption~\eqref{eq-defn} on the measure $e^{-V}$ has been studied in~\cite{BGG11} under the name of $WJ( C )$ inequality; there practical criteria have been given for the inequality to hold. Observe that we can always assume that $C \geq \alpha$ since, if $\alpha \geq 0$, then $\mu$ satisfies a $WJ(\alpha)$ inequality.

\smallskip

\begin{eproof}~
Existence of a stationary solution $\mu_{\infty}$ in $P_2(\rr^n)$ which has a positive density satisfying
$\mu_{\infty}= Z^{-1} e^{-V - W * \mu_{\infty}}$
is given by Proposition~\ref{prop:H}, $iv$; indeed, by~\cite[Cor.~3.11]{BGG11}, assumption~\eqref{eq-defn} on the measure $e^{-V}$ implies the Talagrand  inequality~\eqref{eq:T2} between the Wasserstein distance and the relative entropy (also called $WH$ or $T_2$, see~\cite[Chap. 22]{villani-book1}), with the same $C$.

Then we let $d\mu (x) = e^{-V(x)} \, dx$ and use the convexity assumptions on $V$ and $W$, the bound on $W$ and the sign conditions on $\beta$ and $C - \alpha$ to get, for all $\nu=\nabla\phi\#\mu_\infty$,
\begin{eqnarray*}
&& J_{V,W}(\nu \vert \mu_\infty)\! \! \!
\\
&\ge&\frac{e^{-K}}{Z} \int(\Delta\phi+\Delta\phi^*(\nabla\phi)-2n)\,d\mu
+\int[(\nabla V(\nabla\phi)-\nabla V) \cdot (\nabla\phi-x)-\alpha |\nabla\phi-x|^2]\,d\mu_\infty 
\\
&&+\alpha \int|\nabla\phi - x|^2d\mu_\infty 
+\frac\beta2\iint|\nabla\phi(x)-\nabla\phi (y)-(x-y)|^2 \, d\mu_\infty(x) \, d\mu_\infty(y)\\
&\ge&\!\! \!\! 
\frac{e^{-K}}{Z} \!\! \int(\Delta\phi+\Delta\phi_t^*(\nabla\phi)-2n)\,d\mu
+\frac{e^{-K}}{Z} \!\! \int[(\nabla V(\nabla\phi)-\nabla V) \!\cdot \!(\nabla\phi-x)-\alpha|\nabla\phi-x|^2]\,d\mu \\
&&
+(\alpha + \beta) \int|\nabla\phi-x|^2d\mu_\infty
- \beta \, \Big\vert \int (\nabla \phi - x ) \, d\mu_{\infty} \Big\vert^2\\
&\ge&
(C - \alpha)  \frac{e^{-K}}{Z} \int|\nabla\phi-x|^2d\mu + (\alpha + \beta)\int|\nabla\phi-x|^2 \, d\mu_\infty \\
&\ge&\tilde C \int|\nabla\phi-x|^2 \, d\mu_\infty(x)=\tilde CW_2^2(\nu,\mu_\infty).
\end{eqnarray*}
\end{eproof}
\noindent

\begin{appendix}

\section{Existence of stationary solutions}

The existence of a minimizer of $F$ has been proved by R.~J.~McCann~\cite{mccann-advances} for strictly convex or radially symmetric convex interaction potentials $W$ (and $V=0$). We adapt his classical compactness-lower semicontinuity argument to our diverse cases:

\begin{eprop}\label{prop:H}
The map $F : P_2(\rr^n) \to \rr \cup \{+\infty \}$ defined by~\eqref{eq:F} for absolutely continuous measures and by $+ \infty$ otherwise achieve its minimum in each of the following cases:
\begin{enumerate}
\item $V=0$, $W$ is convex and $W(x) \geq b \vert x \vert^2 - b'$ for $b, b' >0$;

\item $V(x) \geq a \vert x \vert - a'$ and $W(x) \geq b \vert x \vert^2 - b'$ for $a, a', b, b' >0$;

\item $V(x) \geq a \vert x \vert^2 - a'$ and $W(x) \geq b \vert x \vert^2 - b'$ for $b', a, a' >0, b >-a$;

\item $W$ is bounded from below and $e^{-V} \in P_2(\rr^n)$ satisfies a Talagrand transportation inequality 
\begin{equation}\label{eq:T2}
W_2^2(\nu , e^{-V}) \leq \frac{2}{C} \Big( \int \nu \, \log \, \nu \, dx + \int V \, d\nu \Big), \quad \nu \in P_2(\rr^n).
\end{equation}
\end{enumerate}
\end{eprop}

Then, as in~\cite{cmcv-03}, a minimizer $\mu_{\infty}$ of $F$ has a positive density on $\rr^n$ satisfying
$$
\log \, \mu_{\infty} + V + W * \mu_{\infty} = \lambda \in \rr.
$$

\begin{eproof}
First of all, $\inf_{P_2(\rr^n)} F < + \infty$ since $F(\mu) < + \infty$ for $\mu$ the Lebesgue measure on $[0,1]^n$ for instance. Let then $(\mu_p)_p \in P_2(\rr^n)$ be a minimizing sequence, and assume for a while that $\int \vert x \vert^2 d\mu_p$ is bounded. Then $(\mu_p)_p$ is tight, so up to a subsequence admits a limit $\mu_{\infty}$ for the narrow convergence by the Prohorov Theorem. Moreover
$\int \vert x \vert^2 d\mu_{\infty} \leq \liminf_p \int \vert x \vert^2 d\mu_p < + \infty$
so $\mu_{\infty} \in P_2(\rr^n)$. Finally $\mu_{\infty}$ minimizes $F$ on $P_2(\rr^n)$ by lower semicontinuity.

\smallskip

It remains now to bound $\int \vert x \vert^2 d\mu_p$ by $F(\mu_p)$ in each case:

\smallskip

For $i.$, as in~\cite{mccann-advances}, let $\nabla \varphi_p$ transport $\mu_p$ onto $\mu_p(-.)$ and  let ${\bar{\mu}_p} = \frac{I+\nabla \varphi_p}{2} \# \mu_p$ for $I$ the identity map. Now $W$ is convex, so $F$ is displacement convex, so that $F({\bar{\mu}_p}) \leq (F(\mu_p) + F(\mu_p(-.)))/2 = F(\mu_p)$ and $({\bar{\mu}_p})$ is also a minimizing sequence. Moreover
$\int x d{\bar{\mu}_p}=0$ so
$$
\int \vert x \vert^2 d{\bar{\mu}_p}= \frac{1}{2} \iint \vert x-y \vert^2 d{\bar{\mu}_p}(x) d{\bar{\mu}_p}(y) \leq \frac{1}{2} \iint \frac{1}{b} (W(x-y)+b') d{\bar{\mu}_p}(x) d{\bar{\mu}_p}(y) \leq \frac{F({\bar{\mu}_p})}{b}+ \frac{b'}{2b} \cdot
$$

For $ii.$ we observe that 
$$
F(\mu_p) \geq a \int \vert x \vert \, d\mu_p - a' + b \Big[ \int \vert x \vert^2 d\mu_p - \Big\vert \int x \, d\mu_p \Big\vert^2 \Big] -
 \frac{b'}{2}
  \geq a \int \vert x \vert d\mu_p - a'- \frac{b'}{2};
$$
hence $\int \vert x \vert d\mu_p$ is bounded by the second inequality, and then $\int \vert x \vert^2 d\mu_p$ by the first one.

For $iii.$ we  similarly observe, and by discussing on the sign of $b$, that
$$
F(\mu_p) 
\geq 
a \int \vert x \vert^2 d\mu_p - a' + \frac{1}{2} \Big[\iint ( b \vert x-y \vert^2 -b') d\mu_p(x) d\mu_p(y) \Big] 
\geq
(a + \min (b,0)) \int \vert x \vert^2 d\mu_p - a' - \frac{b'}{2} \cdot
 $$
 
 For $iv.$ we notice that
 \begin{eqnarray*}
\int \vert x \vert^2 d\mu_p \leq 2 W_2^2(\mu_p, e^{-V}) + 2  \int \vert x \vert^2 e^{-V} 
&\leq&
 \frac{4}{C} \Big( \int \mu_p \, \log \, \mu_p \, dx + \int V \, d\mu_p \Big) + 2  \int \vert x \vert^2 e^{-V}
\\
&\leq& 
\frac{4}{C} \Big( F(\mu_p) - \frac{1}{2} \inf W \Big) + 2  \int \vert x \vert^2 e^{-V}.
\end{eqnarray*}
\end{eproof}

\end{appendix}

\bigskip

\noindent
{\bf Acknowledgements.} This research was supported  by the French ANR project EVOL.

\newcommand{\etalchar}[1]{$^{#1}$}
\footnotesize{

}


\begin{thebibliography}{10}

\bibitem{ambrosio-gigli-savare}
L.~Ambrosio, N.~Gigli, and G.~Savar{\'e}.
\newblock {\em Gradient flows in metric spaces and in the space of probability
  measures}.
\newblock Lectures in Math. ETH Z\"urich. Birkh\"auser, Basel, 2008.

\bibitem{bbcg08}
D.~Bakry, F.~Barthe, P.~Cattiaux, and A.~Guillin.
\newblock A simple proof of the {P}oincar\'e inequality for a large class of
  probability measures including the log-concave case.
\newblock {\em Elec. Comm. Prob.}, 13:60--66, 2008.

\bibitem{BCCP}
D.~Benedetto, E.~Caglioti, J.~A. Carrillo, and M.~Pulvirenti.
\newblock A non {M}axwellian steady distribution for one-dimensional granular
  media.
\newblock {\em J. Stat. Physics}, 91(5/6):979--990, 1998.

\bibitem{BGG11}
F.~Bolley, I.~Gentil, and A.~Guillin.
\newblock Convergence to equilibrium in {W}asserstein distance for
  {F}okker-{P}lanck equations.
\newblock Preprint, 2011.

\bibitem{bgm}
F.~Bolley, A.~Guillin, and F.~Malrieu.
\newblock Trend to equilibrium and particle approximation for a weakly
  selfconsistent {V}lasov-{F}okker-{P}lanck equation.
\newblock {\em Math. Mod. Num. Anal.}, 44 (5):867--884, 2010.

\bibitem{calvez-carrillo10}
V.~Calvez and J.~A. Carrillo.
\newblock Refined asymptotics for the subcritical {K}eller-{S}egel system and
  related functional inequalities.
\newblock To appear on \emph{Proc. AMS}, 2012.

\bibitem{cmcv-03}
J.~A. Carrillo, R.~J. McCann, and C.~Villani.
\newblock Kinetic equilibration rates for granular media and related equations:
  entropy dissipation and mass transportation estimates.
\newblock {\em Rev. Mat. Iberoam.}, 19(3):971--1018, 2003.

\bibitem{cmcv-06}
J.~A. Carrillo, R.~J. McCann, and C.~Villani.
\newblock Contractions in the 2-{W}asserstein length space and thermalization
  of granular media.
\newblock {\em Arch. Rational Mech. Anal.}, 179:217--263, 2006.

\bibitem{CT03}
J.~A. Carrillo and G.~Toscani.
\newblock Wasserstein metric and large-time asymptotics of nonlinear diffusion
  equations.
\newblock In {\em New trends in math. physics}. World Sci., Singapore, 2005.

\bibitem{portoercole}
J.~A. Carrillo and G.~Toscani.
\newblock Contractive probability metrics and asymptotic behavior of
  dissipative kinetic equations.
\newblock {\em Riv. Mat. Univ. Parma}, 7(6):75--198, 2007.

\bibitem{cgm-08}
P.~Cattiaux, A.~Guillin, and F.~Malrieu.
\newblock Probabilistic approach for granular media equations in the non
  uniformly convex case.
\newblock {\em Prob. Theor. Rel. Fields}, 140(1-2):19--40, 2008.

\bibitem{cordero-gangbo-houdre}
D.~Cordero-Erausquin, W.~Gangbo, and C.~Houdr\'e.
\newblock Inequalities for generalized entropy and optimal transportation.
\newblock In {\em Recent Advances in the Theory and Applications of Mass
  Transport, Contemp. Math. 353}. A. M. S., Providence, 2004.

\bibitem{daneri-savare}
S.~Daneri and G.~Savar\'e.
\newblock Lecture notes on gradient flows and optimal transport.
\newblock To appear on \emph{{S}\'eminaires et {C}ongr\`es, {S}{M}{F}}, 2012.

\bibitem{lisini}
S.~Lisini.
\newblock Nonlinear diffusion equations with variable coefficients as gradient
  flows in {W}asserstein spaces.
\newblock {\em ESAIM Contr. Opt. Calc. Var.}, 15:712--740, 2009.

\bibitem{malrieu03}
F.~Malrieu.
\newblock Convergence to equilibrium for granular media equations and their
  {E}uler schemes.
\newblock {\em Ann. Appl. Probab.}, 13(2):540--560, 2003.

\bibitem{mccann-advances}
R.~J. McCann.
\newblock A convexity principle for interacting gases.
\newblock {\em Adv. Math.}, 128:153--179, 1997.

\bibitem{tugaut2012}
J.~Tugaut.
\newblock Convergence to the equilibria for self-stabilizing processes in
  double well landscape.
\newblock {\em {\rm To appear on} Ann. Prob.}, 2012.

\bibitem{tugaut}
J.~Tugaut.
\newblock Self-stabilizing processes in multi-wells landscape in ${R}^d$ -
  {I}nvariant probabilities.
\newblock {\em {\rm To appear on} J. Theor. Prob.}, 2012.

\bibitem{villani-book1}
C.~Villani.
\newblock {\em Optimal transport, Old and new}, volume 338 of {\em Grund. Math.
  Wiss.}
\newblock Springer, Berlin, 2009.

\end{thebibliography}

\end{document}